\documentclass [a4paper, 12pt]{article}
\usepackage{url}
\usepackage{amsmath}
\usepackage{amssymb}
\usepackage{amsthm}
\usepackage{pgfplots}
\usepackage{changepage}
\usepackage{cite}
\usepackage{geometry}
\pgfplotsset{compat=newest}
\usepackage[colorlinks]{hyperref}
\hypersetup{
	linkcolor = [rgb]{0,0,1},
	urlcolor = [rgb]{0,0,1}
}

\pgfplotsset{every axis/.append style={
		axis x line=middle,    
		axis y line=middle,    
		axis line style={<->}, 
		xlabel={$x$},          
		ylabel={$y$},          
		ticks=none,
}}
\pgfplotsset{mystyle/.style={color=blue,no marks,line width=1pt,<->}} 
\tikzset{>=stealth}

\begin {document}
\title {Finding Elliptic Curves With Many Integral Points}
\author {Benjamin Jones}
\date {\today}
\maketitle

\section {Introduction}

In this paper we outline a method of producing reduced elliptic curves with many integral points and provide the results of the outlined computations, including several curves with hundreds of integral points. The first three sections give background and describe our work with integral points on elliptic curves. The last section is unrelated to elliptic curves and provides a complete classification of self-descriptive numbers.

\subsection {What is an elliptic curve?}

The brief overview of elliptic curves we give here should be sufficient for understanding and motivating the results we give later. For the purposes of this paper, an \emph{elliptic curve} $E$ is a curve on the rational plane of the form:
$$E: Y^2 + aXY + cY = X^3 + bX^2 + dX + e$$
where $a$, $b$, $c$, $d$, and $e$ are integers. A \emph{point} on the elliptic curve is a rational pair $(x, y)$ which satisfies the equation. If $(x, y)$ is a point on the curve and $x$ and $y$ are integers, then the point is called an \emph{integral point} on the elliptic curve.

Suppose we have an elliptic curve $E$:
$$E: Y^2 + aXY + cY = X^3 + bX^2 + dX + e$$
and some integer $g$ where $e = g^6e'$, $d = g^4d'$, $c = g^3c'$, $b = g^2b'$, and $a = ga'$. Then if $(x, y)$ is an integral point on $E$, the point $(\frac{x}{g^2}, \frac{y}{g^3})$ lies on the curve:
$$g^6Y^2 + g^6a'XY + g^6c'Y = g^6X^3 + g^6b'X^2 + g^6d'X + g^6e'$$
and hence, $(\frac{x}{g^2}, \frac{y}{g^3})$ is a point on the elliptic curve $E'$:
$$E': Y^2 + a'XY + c'Y = X^3 + b'X^2 + d'X + e'$$
We have a bijection between the points on $E$ and the points on $E'$ given by:
$$(x, y) \leftrightarrow (\frac{x}{g^2}, \frac{y}{g^3})$$
In general, the more factors $g$ has, the higher the number of integral points on $E$.

If for any elliptic curve $E$ such a $g > 1$ exists, then we call $E$ an \emph{unreduced} curve. Otherwise, we call $E$ \emph{reduced}. By creating an unreduced elliptic curve $E$ which reduces to a reduced curve $E'$, we can artificially increase the number of integral points on $E$ by increasing $g$. Therefore in this paper we are only concerned with finding reduced elliptic curves.

\subsection {Groups}

We will now define a group. Let $G$ be any set, and $*$ be any binary operation on $G$. Then $(G, *)$ is a group if the following hold:

\begin{enumerate}
	\item There is some $e \in G$ such that for all $g \in G$, $e*g = g*e = g$. We call $e$ the identity of the group.
	\item For every $g \in G$ there exists some inverse element denoted $g^{-1} \in G$ such that $g*g^{-1} = g^{-1}*g = e$.
	\item For every $g_1, g_2, g_3 \in G$ we have $(g_1*g_2)*g_3 = g_1*(g_2*g_3)$. In other words, $*$ is associative.
\end{enumerate}

Often the group $(G, *)$ is abbreviated as $G$ when the binary operation is understood. A group is called \emph{abelian} if for every $a \in G$ and $b \in G$, $a*b = b*a$. When $G$ is abelian, we often use $a + b$ to denote $a*b$, and when $G$ is not abelian, $ab$ denotes $a*b$. The \emph{order} of $G$ is the cardinality of the set $G$. If $g \in G$ and $n \in \mathbb{N}$ with $n > 0$, then $ng$ when $G$ is abelian, or $g^n$ when $G$ is not abelian, denotes $\underbrace{g * g * g * ... * g}_\text{n times}$. If there is an $n$ such that $ng = e$ or $g^n = e$, then the smallest such $n$ is called the \emph{order} of $g$, and g has \emph{finite order}. If no such $n$ exists then g has \emph{infinite order}.

A set $S \subseteq G$ \emph{generates} $G$ if every $g \in G$ can be written as $s_1 * s_2 * s_3 * ... * s_n$ where for each $i$, either $s_i \in S$ or $s_i^{-1} \in S$. The elements of $S$ are \emph{generators} of $G$ if $S$ generates $G$ and no proper subset of $S$ generates $G$. A group $G$ is \emph{finitely generated} if it can be generated by some finite subset $S$ of $G$.

\subsection {Group Structure on Elliptic Curves}

In this section we will define a group structure on the rational points of non-singular elliptic curves. To see how, let's suppose we have two points on a non-singular elliptic curve $E$:
$$E: Y^2 + aXY + cY = X^3 + bX^2 + dX + e$$
$$P = (x_P, y_P), Q = (x_Q, y_Q)$$

For now, suppose $P$ and $Q$ are distinct points on $E$. If $x_P \not = x_Q$ then the line $L$ which passes through both $P$ and $Q$ can be written as:
$$L: Y = m(X - x_P) + y_P$$
where $m = \frac{y_P - y_Q}{x_P - x_Q}$ is rational. Then, if $(X, Y)$ lies on both $E$ and $L$, we can substitute this expression for $Y$ to get:
$$X^3 + bX^2 + dX + e - [m(X - x_P + y_P)]^2 - aX[m(X - x_P) + y_P] - c[m(X - x_P) + y_P] = 0$$
Since the left-hand side is a cubic in $X$, it has at most $3$ roots. Also, since $P$ and $Q$ lie on $E$ and $L$, we know that $x_P$ and $x_Q$ are two rational roots of the cubic. Thus the third root of the cubic, which we will call $x_R$, must be rational. Once we compute $x_R$, we can use the equation for $L$ to compute the point $R = (x_R, y_R)$ which lies on both $L$ and $E$. If $x_R = x_P$ or $x_R = x_Q$, then $R = P$ or $R = Q$, and the line $L$ is tangent to $E$. Otherwise, we see that the line between $P$ and $Q$ crosses the elliptic curve $E$ at exactly one other rational point.

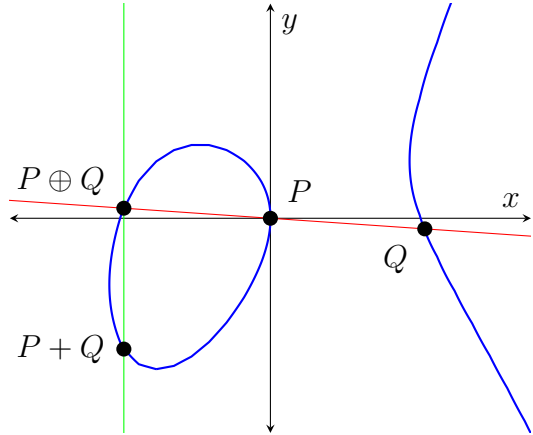
\begin{figure}
\centering
\begin{tikzpicture}
	\begin{axis}[
		xmin=-3,xmax=3,
		ymin=-3,ymax=3,
	]
	\addplot +[no markers,
		raw gnuplot,
		thick,
		empty line = jump,
		blue
		] gnuplot {
		set contour base;
		set cntrparam levels discrete 0.003;
		unset surface;
		set view map;
		set isosamples 500;
		splot y^2 - x*y - x^3 + 3*x;
		};
	\addplot +[mark = none,red] plot {x*-0.0833333333};
	\addplot +[mark = node,green] plot ({-27/16}, {x});
	\node [label={45:{$P$}}, circle, fill, inner sep=2pt] at (axis cs:0, 0) {};
	\node [label={225:{$Q$}}, circle, fill, inner sep=2pt] at (axis cs:16/9, -4/27) {};
	\node [label={165:{$P\oplus Q$}}, circle, fill, inner sep=2pt] at (axis cs: -27/16, 9/64) {};
	\node [label={180:{$P + Q$}}, circle, fill, inner sep=2pt] at (axis cs: -27/16, -117/64) {};
	\end{axis}
\end{tikzpicture}
\caption{When $P \not = Q$ and $x_P \not = x_Q$}
\end{figure}

If $x_P = x_Q$, however, then this is not true. In this case, the intersection between the line crossing through $P$ and $Q$ and the curve $E$ consists of all points $(x_P, Y)$ where $Y$ satisfies the equation:
$$Y^2 + ax_PY + cY - x_P^3 - bx_P^2 - dx_P - e = 0$$
This is a quadratic in $Y$, and can only have at most two solutions. However, $y_P$ and $y_Q$ are already two distinct solutions to this quadratic, so there must be no third distinct intersection point.

So when $x_P \not = x_Q$ we will define $P \oplus Q$ as the rational point $R$ created from the intersection between the line from $P$ and $Q$ and the curve $E$. For any point $P$ in $E(\mathbb{Q})$, we will define $-P$ as the other intersection between the vertical line passing through $P$ and the curve $E$ if it exists, otherwise $-P = P$. Then we define $P + Q$ to be $-(P \oplus Q)$. In other words, if $P + Q = (x', y')$ then $x' = x_R$ and $y'$ is other root of this quadratic in $Y$:
$$Y^2 + ax_RY + cY - x_R^3 - bx_R^2 - dx_R - e$$

\begin{figure}
\centering
\begin{tikzpicture}
	\begin{axis}[
		xmin=-2,xmax=4,
		ymin=-5,ymax=4,
	]
	\addplot +[no markers,
		raw gnuplot,
		thick,
		empty line = jump,
		blue
		] gnuplot {
		set contour base;
		set cntrparam levels discrete 0.003;
		unset surface;
		set view map;
		set isosamples 500;
		splot (y^2 - x*y - x^3 + 3*x;
		};
	\addplot +[mark = none,red] plot {(x - 3)*-2.3333333333 - 3};
	\addplot +[mark = node,green] plot ({16/9}, {x});
	\node [label={180:{$P$}}, circle, fill, inner sep=2pt] at (axis cs:3, -3) {};
	\node [label={195:{$P\oplus P$}}, circle, fill, inner sep=2pt] at (axis cs: 16/9, -4/27) {};
	\node [label={0:{$2P$}}, circle, fill, inner sep=2pt] at (axis cs: 16/9, 52/27) {};
	\end{axis}
\end{tikzpicture}
\caption{When $P = Q$}
\end{figure}
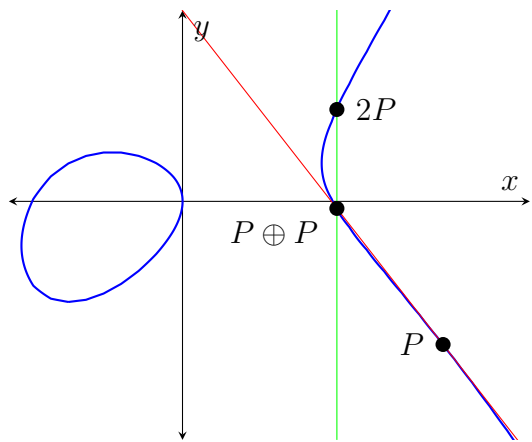

When $x_P = x_Q$, we have $P = -Q$ so $P + Q$ is the identity. However, the line passing through $P$ and $Q$ doesn't intersect the curve at a third point. Because of this, for the identity we introduce a new `point at infinity' denoted $\mathcal{O}$. In projective coordinates, we could write the curve $E$ as:
$$Y^2Z + aXYZ + cYZ^2 = X^3 + bX^2Z + dXZ^2 + eZ^3$$
Then the point $\mathcal{O}$ would represent the coordinate with $Z = 0$. Alternatively, one can add $\mathcal{O}$ as a formal symbol to the set of rational points of $E$ on the plane. We of course define $-\mathcal{O} = \mathcal{O}$ and $\mathcal{O} + S = S + \mathcal{O} = S$ for any rational point $S$.

Finally, if $P = Q$ we instead choose $L$ to be the line tangent to $E$ at the $P$. Since $E$ is non-singular, an unambiguous choice for the slope of the tangent line at $P$ does exist. Again, if this line is not vertical, then it intersects $E$ at one additional point. We define this other point to be $P \oplus P$, and $P + P = 2P$ to be $-(P \oplus P)$.

The group of rational points on a non-singular elliptic curve $E$ along with the point $\mathcal{O}$ will be denoted $E(\mathbb{Q})$. All of the axioms of a group follow quickly from the above description except associativity. For associativity, refer to \cite[pg.19]{silvermanec} and \cite[pg.28]{casselsec}.

\subsection {Important Theorems}

\emph{Mordell's Theorem}: If $E$ is a non-singular elliptic curve, then there is a finite set of generators for $E(\mathbb{Q})$. \cite[pg.54]{casselsec} \cite[pg.63]{silvermanec}

If $S$ is a minimal set of generators for $E(\mathbb{Q})$, then $S$ is a finite set and we call the number of generators in this set with infinite order the \emph{rank} of $E$. In particular, the rank of $E$ is always finite. If the rank of $E(\mathbb{Q})$ is $0$, then $E$ has only finitely many rational points. Otherwise, $E$ has an infinite number of rational points.

\emph{Siegel's Theorem}: If $E$ is a non-singular elliptic curve, then $E$ has only finitely many integral points. \cite[pg.146]{silvermanec}

\section{The Parameterization}
First we may start with an elliptic curve in Weirstrass form:
$$y^2 + axy + cy = x^3 + bx^2 + dx + f$$
where $a$, $b$, $c$, $d$, and $f$ are integers.
For convenience, we can center the curve over the point $(0, 0)$, in which case we would have $f = 0$ and the curve:
$$y^2 + axy + cy = x^3 + bx^2 + dx$$
Let $P$ refer to the point $(0, 0)$ on this curve.

\subsection{2P Integral}

Now we want to impose the restriction that $2P$ is integral. On curves of this form, $2(x, y) = (x', y')$ is given by

\begin{equation}
\lambda = \frac {3x^2 + 2bx - ay + d} {2y + ax + c}
\end{equation}
\begin{equation}
x' = \lambda^2 + \lambda a - b - 2x
\end{equation}
\begin{equation}
y' = -ax' - c - \lambda x' + \lambda x_1 - y
\end{equation}

In the case for $(x, y) = P = (0, 0)$ we have:

$$\lambda = \frac {d} {c}$$

$$x' = \frac {d^2 + acd} {c^2} - b$$

If $c | d$ then $x'$ is an integer and hence $(x', y')$ is integral. On the other hand, if $x'$ is an integer then we would have that $\frac {d^2 + acd} {c^2} = (\frac {d} {c})^2 + a(\frac {d} {c})$ is an integer. We want to show that $d | c$. More generally, we need to prove the following:
\\\\
\emph{Lemma}. Suppose $q = \frac {m} {n}$ is rational with $gcd(m, n) = 1$ and that $w$ is an integer. Then $q^2 + wq$ is an integer if and only if $q$ is an integer.
\begin{proof}
It is clear that if $q$ is an integer then $q^2 + wq$ is an integer. Now suppose $q^2 + wq$ is an integer. Then
$$q^2 + wq = \left(\frac {m} {n}\right)^2 + w\cdot \frac{m} {n} = \frac {m^2 + wmn} {n^2} = m\cdot \frac {m + wn} {n^2}$$
Is an integer. This is only true if $n^2 | m + wn$. Since $n | wn$, we must have that $n | m$. Since $gcd(n, m) = 1$, we have that $n = 1$, and $q$ is an integer.
\end{proof}
This lemma will be used frequently in the parameterization up to $5P$. By the lemma, if $x'$ is an integer then $\frac{c} {d}$ is an integer and $d | c$. So $2P$ is integral if and only if $d | c$. Let $e = \frac {d} {c}$. Then our curve has the form:
$$y^2 + axy + cy = x^3 + bx^2 + ecx$$
This is a parameterization of all curves with $2P$ integral. As it happens, there is a transformation which will simplify the above. After substituting $ex + y$ for $y$ and simplifying, we obtain the equivalent curve:
$$y^2 + (2e + a)xy + cy = x^3 + (b - e^2 - ae)x^2$$
Since the transformation did not change the x coordinate, a point in this curve is integral if and only if it is integral in the previous curve. In addition, this curve is clearly aleady in the form of the previous curve. Therefore, every curve for which $2P$ is integral can be transformed into an equivalent curve for which $e = 0$. Without loss of generality, we can assume $e = 0$, and obtain the parameterization:
$$y^2 + axy + cy = x^3 + bx^2$$
for curves which have $2P$ integral. On any such curve

$$2P = (-b, ab - c)$$

\subsection{3P Integral}

On any curve of this form, we can compute $(x', y') = (x, y) + (0, 0)$ as:

$$x' = \left(\frac {y} {x}\right)^2 + a\left(\frac {y} {x}\right) - b$$
$$y' = -x'\left(a + \frac {y} {x}\right) - c$$

Since $b$ is an integer, our previous lemma tells us that $(x', y')$ is integral if and only if $x | y$. Thus, $3P = (-b, ab - c) + (0, 0)$ is integral if and only if $b | ab - c$, which is true if and only if $b | c$. Replacing $c$ with $db$, we obtain the parameterization:

$$y^2 + axy + dby = x^3 + bx^2$$

for all curves with $3P$ integral. For this curve

$$3P = (d^2 - ad, -d^3 + ad^2 - bd)$$

\subsection{4P Integral}

Next, we can compute $4P = (d^2 - ad, -d^3 + ad^2 - bd) + (0, 0)$. The $x$ coordinate is equal to:

$$\left(\frac{b} {a - d}\right)^2 - d\cdot \frac{b} {a - d}$$

By our lemma, this is an integer if and only if $a - d | b$. So we can substitute $f(a - d)$ for $b$ in the parameterization to get:

$$y^2 + axy + df(a - d)y = x^3 + f(a - d)x^2$$

which parameterizes the curves for which $4P$ is integral. On this curve, $4P = (f^2 - df, -f^3 + (2d - a)f^2)$.

\subsection {5P Integral}

The $x$ coordinate of $5P$ is:

$$\left(\frac {d(a - d)} {d - f}\right)^2 + a\left(\frac {d(a - d)} {d - f}\right) + d(a - d)$$
Which is integral if and only if $d - f|d(a - d)$. Suppose $d(a - d) = g(d - f)$ for some integer $g$. Then, let $t = gcd(d, g)$. For some $r$ and $u$, we have $g = tu$ and $d = tr$ where $gcd(u, r) = 1$. Substituting $tu$ and $tr$ for $g$ and $d$, we have that:
$$tru - fu = ra - r^2t$$
Which implies that $r | fu$. But since $gcd(r, u) = 1$, we must have that $r | f$. So write $f = rs$ for some integer $s$. Substituting $rs$ for $f$ yields:
$$(tr - rs)tu = tr(a - tr)$$
Solving for $a$ gives:
$$a = tr + tu - su$$
Substituting for $d$ and $f$ as well in the previous parameterization gives us a new parameterization:
$$y^2 + (tr + tu - su)xy + tr^2su(t - s)y = x^3 + rsu(t - s)x^2$$
and on this curve, $5P$ is the point:
$$(stu^2 - rstu, rs^2tu^2 - st^2u^3)$$
This is the last parameterization we will give which is general. Any elliptic curve of the form $$y^2 + axy + cy = x^3 + bx^2 + dx + e$$ with $2(0,0)$, $3(0,0)$, $4(0,0)$, and $5(0,0)$ integral can be reduced to an equivalent curve in the form of our parameterization. Notice that if $g = gcd(t, s) \ne 1$ then the curve given by the parameterization is unreduced, so we only need to consider parameters for which $gcd(t, s) = 1$.

\subsection {8P Integral}

From the previous parameterization, it is possible to derive a parameterization for curves in which $2P$ up to $8P$ are integral. However, it is not a parameterization of all such curves as are the previous ones. The parameterization has the following form:
$$y^2 + axy + cy = x^3 + bx^2$$
with:
$$a = u^2 + puv + u^2v - p^2v^2 + 2puv^2 - puv^3 + u^2v^3$$
$$b = uv(v + 1)(u - p)(u + pv)(u + uv + pv + uv^2)$$
$$c = bv^2(u - p)(p + uv)$$
One can verify using \emph{Magma} that the $x$ coordinate of $8P$ on this curve is:
$$-u^3v^5p - u^3v^4p - u^3v^3p + u^2v^5p^2 + u^2v^3p^2 + uv^4p^3$$
And $8P$ is integral.

\section{Search Results}
\subsection{Highest Multiple}

One may ask what the highest integral multiple of an integral point on a reduced curve is. Searching through all reduced elliptic curves to find such a point and curve would be inefficient. Our parameterization of curves for which $nP$ is integral guarantees that the largest integral multiple of $P$ is at least $nP$. These parameterizations therefore are a better basis for a computer search.

The search was done as follows:
\begin{enumerate}
\item A C program using GMP searched through curves in the parameterization for $5P$ integral. For each curve it checked whether $nP$ was integral for $1 \leq n \leq 35$. If for some $n \geq 15$, $nP$ was integral, the curve was added as a candidate.
\item A larger search using GMP over the parameterization for $8P$ integral was performed, and it used the same criterion to add curves as a candidate.
\item In \emph{Magma}, the curves were checked against their minimal model and sorted based on the highest integral multiple of $P$.
\end{enumerate}

Implementing the search in C instead of \emph{Magma} made the search much faster, and checking candidate curves in \emph{Magma} allowed us to use its more robust tools for sorting and verification of candidate curves. Below are the curves we found with the highest multiples of $P$:

\begin{adjustwidth}{-0.75in}{-0.75in}
\begin{tabular}{l|l}
Curve & Integral Multiples of (0,0) \\
$y^2 - 17xy + 960y = x^3 - 30x^2                 $&$ 1, 2, 3, 4, 5, 6, 7, 8, 9, 10, 11, 14, 15, 17, 31      $\\
$y^2 + 7xy - 210y = x^3 + 70x^2                  $&$ 1, 2, 3, 4, 5, 6, 7, 9, 10, 11, 12, 14, 15, 18, 25     $\\
$y^2 - 28xy - 840y = x^3 - 420x^2                $&$ 1, 2, 3, 4, 6, 7, 8, 9, 11, 12, 13, 24                 $\\
$y^2 + 535xy + 10929600y = x^3 + 22770x^2        $&$ 1, 2, 3, 4, 5, 6, 8, 9, 10, 11, 12, 15, 23, 24         $\\
$y^2 + 1879xy - 300699000y = x^3 - 155400x^2     $&$ 1, 2, 3, 4, 5, 6, 8, 9, 10, 11, 12, 22                 $\\
$y^2 - 80xy - 34560y = x^3 + 480x^2              $&$ 1, 2, 3, 4, 5, 6, 7, 9, 10, 11, 12, 21                 $\\
$y^2 - 77xy - 7920y = x^3 - 2640x^2              $&$ 1, 2, 3, 4, 5, 6, 7, 8, 9, 12, 13, 14, 15, 16, 18, 21  $\\
$y^2 + 1107xy + 102316500y = x^3 + 104940x^2     $&$ 1, 2, 3, 4, 6, 7, 8, 9, 12, 21                         $\\
$y^2 - 181xy - 436800y = x^3 + 2730x^2           $&$ 1, 2, 3, 4, 5, 6, 7, 8, 9, 10, 11, 12, 15, 21          $\\
$y^2 + 253xy + 1197000y = x^3 + 5320x^2          $&$ 1, 2, 3, 4, 5, 6, 7, 9, 10, 12, 14, 15, 18, 19, 20, 21 $\\
$y^2 + 211xy + 537030y = x^3 + 6630x^2           $&$ 1, 2, 3, 4, 5, 6, 7, 8, 9, 11, 12, 13, 15, 16, 21      $\\
$y^2 - 3599xy - 116180064y = x^3 - 4149288x^2    $&$ 1, 2, 3, 4, 5, 6, 7, 8, 10, 21                         $\\
$y^2 + 11xy + 1050y = x^3 + 210x^2               $&$ 1, 2, 3, 4, 5, 6, 7, 8, 9, 10, 12, 13, 15, 17, 21      $\\
$y^2 - 2813xy - 2968105140y = x^3 + 19399380x^2  $&$ 1, 2, 3, 4, 5, 6, 7, 8, 9, 10, 12, 21                  $\\
$y^2 + 479xy + 12061500y = x^3 + 43860x^2        $&$ 1, 2, 3, 4, 5, 6, 7, 8, 9, 12, 21                      $\\
$y^2 + 1543xy + 262765440y = x^3 + 191520x^2     $&$ 1, 2, 3, 4, 5, 6, 7, 8, 9, 11, 21                      $\\
$y^2 + 99xy + 1928934y = x^3 - 30618x^2          $&$ 1, 2, 3, 4, 5, 6, 7, 8, 10, 11, 12, 14, 15, 20         $\\
$y^2 - 3xy - 12096y = x^3 - 672x^2               $&$ 1, 2, 3, 4, 5, 6, 8, 9, 10, 11, 20                     $\\
$y^2 - 133xy - 49686y = x^3 + 546x^2             $&$ 1, 2, 3, 4, 5, 6, 7, 8, 10, 12, 13, 15, 16, 19, 20     $\\
$y^2 - 3659xy - 51856042050y = x^3 - 27422550x^2 $&$ 1, 2, 3, 4, 5, 6, 7, 8, 10, 12, 13, 20                 $\\
$y^2 - 133xy - 35190y = x^3 - 7038x^2            $&$ 1, 2, 3, 4, 5, 6, 9, 10, 11, 15, 20                    $\\
$y^2 + 2921xy - 98463750y = x^3 - 112530x^2      $&$ 1, 2, 3, 4, 5, 6, 7, 9, 10, 11, 14, 20                 $\\
$y^2 - 3141xy - 1765062090y = x^3 + 614790x^2    $&$ 1, 2, 3, 4, 5, 6, 8, 9, 10, 12, 15, 18, 20             $\\
$y^2 - 65xy - 1518y = x^3 - 1518x^2              $&$ 1, 2, 3, 4, 5, 6, 9, 10, 12, 15, 18, 20                $\\
$y^2 + 707xy + 27550320y = x^3 + 49910x^2        $&$ 1, 2, 3, 4, 5, 7, 8, 10, 20                            $\\
$y^2 + 718xy - 347490y = x^3 - 115830x^2         $&$ 1, 2, 3, 4, 5, 8, 9, 10, 11, 20                        $\\
$y^2 + 123xy + 50540y = x^3 + 532x^2             $&$ 1, 2, 3, 4, 5, 7, 8, 10, 20                            $\\
$y^2 - 151xy - 184800y = x^3 - 13200x^2          $&$ 1, 2, 3, 4, 5, 6, 7, 8, 9, 10, 11, 20                  $\\
$y^2 + 5273xy + 12331371150y = x^3 + 2920050x^2  $&$ 1, 2, 3, 4, 5, 7, 8, 10, 20                            $\\
$y^2 + 341xy + 2827440y = x^3 + 9240x^2          $&$ 1, 2, 3, 4, 5, 6, 8, 9, 10, 12, 20                     $\\
$y^2 + 103xy + 26730y = x^3 + 330x^2             $&$ 1, 2, 3, 4, 5, 6, 7, 9, 10, 11, 13, 15, 20             $\\
$y^2 - 157xy + 288990y = x^3 - 1170x^2           $&$ 1, 2, 3, 4, 5, 6, 7, 8, 10, 12, 13, 14, 17, 20         $\\
$y^2 - 5xy + 48y = x^3 - 6x^2                    $&$ 1, 2, 3, 4, 5, 6, 7, 8, 11, 19                         $\\
$y^2 - 396xy - 3732480y = x^3 + 10368x^2         $&$ 1, 2, 3, 4, 5, 6, 8, 9, 10, 19                         $\\
$y^2 + 35xy + 1650y = x^3 + 330x^2               $&$ 1, 2, 3, 4, 5, 6, 7, 8, 9, 12, 19                      $\\
$y^2 + 3739xy - 777288960y = x^3 + 456960x^2     $&$ 1, 2, 3, 4, 5, 6, 7, 8, 11, 19                         $\\
$y^2 + 213xy - 5738040y = x^3 - 20790x^2         $&$ 1, 2, 3, 4, 6, 8, 9, 18                                $\\
$y^2 + 1981xy - 121307340y = x^3 - 1002540x^2    $&$ 1, 2, 3, 4, 5, 6, 7, 8, 9, 12, 18                      $\\
$y^2 + 161xy - 36960y = x^3 + 9240x^2            $&$ 1, 2, 3, 4, 5, 6, 7, 8, 9, 11, 12, 18                  $\\
$y^2 - 24xy - 1634239152y = x^3 + 1464372x^2     $&$ 1, 2, 3, 6, 9, 18                                     $
\end{tabular}
\end{adjustwidth}

\subsection{Most Multiples Integral}

Not only are the parameterizations a natural way of finding curves with high integral multiples, but they also make it easier to find curves with many integral multiples of a point. For instance, the parameterization for $8P$ guarantees that the curves have at least $8$ integral multiples of $P$. Candidate curves produced for the last search were sorted with regards to the number of integral multiples of $P$. Below are the curves with the highest number of integral multiples of $P$:

\newgeometry{margin=1in}
\begin{adjustwidth}{-0.75in}{-0.75in}
\begin{tabular}{p{0.6\linewidth} | p{0.4\linewidth}}
Curve & Integral Multiples of (0, 0) \\
$y^2 + 253xy + 1197000y = x^3 + 5320x^2                                                   $&$ 1, 2, 3, 4, 5, 6, 7, 9, 10, 12, 14, 15, 18, 19, 20, 21 $\\
$y^2 - 77xy - 7920y = x^3 - 2640x^2                                                       $&$ 1, 2, 3, 4, 5, 6, 7, 8, 9, 12, 13, 14, 15, 16, 18, 21  $\\
$y^2 - 17xy + 960y = x^3 - 30x^2                                                          $&$ 1, 2, 3, 4, 5, 6, 7, 8, 9, 10, 11, 14, 15, 17, 31      $\\
$y^2 + 11xy + 1050y = x^3 + 210x^2                                                        $&$ 1, 2, 3, 4, 5, 6, 7, 8, 9, 10, 12, 13, 15, 17, 21      $\\
$y^2 + 211xy + 537030y = x^3 + 6630x^2                                                    $&$ 1, 2, 3, 4, 5, 6, 7, 8, 9, 11, 12, 13, 15, 16, 21      $\\
$y^2 + 1087xy - 2063880y = x^3 - 294840x^2                                                $&$ 1, 2, 3, 4, 5, 6, 7, 8, 9, 10, 11, 12, 15, 16, 18      $\\
$y^2 + 7xy - 210y = x^3 + 70x^2                                                           $&$ 1, 2, 3, 4, 5, 6, 7, 9, 10, 11, 12, 14, 15, 18, 25     $\\
$y^2 - 209xy - 2446080y = x^3 + 23520x^2                                                  $&$ 1, 2, 3, 4, 5, 6, 7, 8, 9, 10, 11, 12, 13, 14, 18      $\\
$y^2 - 133xy - 49686y = x^3 + 546x^2                                                      $&$ 1, 2, 3, 4, 5, 6, 7, 8, 10, 12, 13, 15, 16, 19, 20     $\\
$y^2 - 157xy + 288990y = x^3 - 1170x^2                                                    $&$ 1, 2, 3, 4, 5, 6, 7, 8, 10, 12, 13, 14, 17, 20         $\\
$y^2 + 863xy + 47278080y = x^3 + 61560x^2                                                 $&$ 1, 2, 3, 4, 5, 6, 7, 8, 9, 11, 12, 14, 15, 16          $\\
$y^2 + 1601xy - 72292500y = x^3 - 64260x^2                                                $&$ 1, 2, 3, 4, 5, 6, 7, 8, 9, 10, 11, 12, 13, 16          $\\
$y^2 - 181xy - 436800y = x^3 + 2730x^2                                                    $&$ 1, 2, 3, 4, 5, 6, 7, 8, 9, 10, 11, 12, 15, 21          $\\
$y^2 - 41xy - 2211300y = x^3 + 5460x^2                                                    $&$ 1, 2, 3, 4, 5, 6, 7, 8, 9, 10, 11, 12, 13, 16          $\\
$y^2 + 73xy + 28600y = x^3 + 440x^2                                                       $&$ 1, 2, 3, 4, 5, 6, 8, 9, 10, 12, 15, 16, 17, 18         $\\
$y^2 - 359xy + 2664750y = x^3 - 6270x^2                                                   $&$ 1, 2, 3, 4, 5, 6, 7, 8, 9, 10, 11, 12, 14, 15          $\\
$y^2 + 99xy + 1928934y = x^3 - 30618x^2                                                   $&$ 1, 2, 3, 4, 5, 6, 7, 8, 10, 11, 12, 14, 15, 20         $\\
$y^2 - 103xy - 2522520y = x^3 - 32760x^2                                                  $&$ 1, 2, 3, 4, 5, 6, 7, 8, 9, 10, 12, 13, 14, 15          $\\
$y^2 - 821xy - 2111469360y = x^3 + 510510x^2                                              $&$ 1, 2, 3, 4, 5, 6, 7, 8, 10, 12, 13, 14, 15, 16         $\\
$y^2 + 535xy + 10929600y = x^3 + 22770x^2                                                 $&$ 1, 2, 3, 4, 5, 6, 8, 9, 10, 11, 12, 15, 23, 24         $\\
$y^2 - 3823xy - 20401123050y = x^3 - 18102150x^2                                          $&$ 1, 2, 3, 4, 5, 6, 7, 8, 9, 10, 11, 12, 14, 15          $\\
$y^2 + 43xy - 2337720y = x^3 + 30360x^2                                                   $&$ 1, 2, 3, 4, 5, 6, 7, 8, 9, 10, 11, 12, 15              $\\
$y^2 + 703xy - 34471710y = x^3 + 166530x^2                                                $&$ 1, 2, 3, 4, 5, 6, 7, 8, 9, 10, 12, 15, 16              $\\
$y^2 + 1583xy + 296229420y = x^3 + 211140x^2                                              $&$ 1, 2, 3, 4, 5, 6, 7, 8, 10, 12, 13, 14, 15             $\\
$y^2 - 1525xy - 6727687200y = x^3 - 6468930x^2                                            $&$ 1, 2, 3, 4, 5, 6, 7, 8, 9, 12, 13, 14, 15              $\\
$y^2 + 53xy - 32340y = x^3 + 4620x^2                                                      $&$ 1, 2, 3, 4, 5, 6, 7, 8, 9, 10, 11, 12, 15              $\\
$y^2 - 211xy - 9077250y = x^3 + 27930x^2                                                  $&$ 1, 2, 3, 4, 5, 6, 7, 9, 10, 11, 12, 14, 15             $\\
$y^2 + 323790750569xy - 64392933999375238312586416005120y = x^3 -198500546018619925080x^2 $&$ 1, 2, 3, 4, 5, 6, 7, 8, 9, 10, 12, 14, 15              $\\
$y^2 + 91xy - 7362630000y = x^3 - 4395600x^2                                              $&$ 1, 2, 3, 4, 5, 6, 7, 8, 10, 12, 13, 14, 15             $\\
$y^2 + 38576xy + 6225494016000y = x^3 + 187514880x^2                                      $&$ 1, 2, 3, 4, 5, 6, 7, 8, 9, 10, 12, 14, 15              $\\
$y^2 + 3193xy + 1880615880y = x^3 + 703560x^2                                             $&$ 1, 2, 3, 4, 5, 6, 7, 8, 9, 10, 11, 12, 15              $\\
$y^2 - 6525xy + 199905468750y = x^3 - 22781250x^2                                         $&$ 1, 2, 3, 4, 5, 6, 7, 8, 10, 11, 12, 14, 17             $\\
$y^2 - 156xy - 268800y = x^3 + 1920x^2                                                    $&$ 1, 2, 3, 4, 5, 6, 7, 8, 10, 12, 13, 14, 15             $\\
$y^2 + 116xy - 5376000y = x^3 - 26880x^2                                                  $&$ 1, 2, 3, 4, 5, 6, 7, 8, 9, 10, 11, 12, 16              $\\
$y^2 - 3141xy - 1765062090y = x^3 + 614790x^2                                             $&$ 1, 2, 3, 4, 5, 6, 8, 9, 10, 12, 15, 18, 20             $\\
$y^2 - 53xy - 10080y = x^3 + 210x^2                                                       $&$ 1, 2, 3, 4, 5, 6, 7, 8, 9, 12, 13, 14, 15              $\\
$y^2 + 103xy + 26730y = x^3 + 330x^2                                                      $&$ 1, 2, 3, 4, 5, 6, 7, 9, 10, 11, 13, 15, 20             $\\
$y^2 - 653xy - 1347570y = x^3 + 58590x^2                                                  $&$ 1, 2, 3, 4, 5, 6, 7, 8, 9, 10, 12, 14, 15              $\\
$y^2 + 93160101824xy + 3534126808484560635939394682880y = x^3 +38102459540853227520x^2    $&$ 1, 2, 3, 4, 5, 6, 7, 8, 9, 10, 12, 14, 15              $\\
$y^2 + 311xy - 42020160y = x^3 - 101010x^2                                                $&$ 1, 2, 3, 4, 5, 6, 7, 8, 9, 10, 12, 14, 15             $
\end{tabular}
\end{adjustwidth}
\clearpage
\restoregeometry

Note that since candidate curves are curves for which we have $n \geq 15$ with $nP$ integral, a curve with only $P$ through $14P$ integral will not appear in the table above. For $n < 15$, there may be curves given by the parameterizations for which $n$ positive multiples of $P$ are integral. However, the first curve in the above table contains $16$ integral points, so this search procedure certainly still yields the best curves given by the parameterizations in the search range. 

\subsection{Smallest Height}

Candidate curves from the previous searches also tend to have a small height. For a definition of the canonical height, see \cite[pg.82]{casselsec} or \cite[pg.103]{silvermanec}. Roughly, the height of a point $P$ is an estimate of how large the denominator of the x and y coordinates of $P$ become with increasing multiples of $P$. Integral multiples of $P$ are multiples for which this denominator is $1$. Since the candidate curves from previous searches have large multiples of $P$ integral, the height of $P$ on these curves is small. Below are the candidate curves sorted by the smallest height of $P$:

\begin{adjustwidth}{-0.75in}{-0.75in}
\begin{tabular}{l|l}
Curve & Height of (0,0) on the Minimal Model \\
$y^2 + 253xy + 1197000y = x^3 + 5320x^2           $&$ 0.008914$ \\
$y^2 - 77xy - 7920y = x^3 - 2640x^2               $&$ 0.009039$ \\
$y^2 + 73xy + 28600y = x^3 + 440x^2               $&$ 0.009740$ \\
$y^2 + 7xy - 210y = x^3 + 70x^2                   $&$ 0.009964$ \\
$y^2 - 8xy - 140y = x^3 - 70x^2                   $&$ 0.011278$ \\
$y^2 - 133xy - 49686y = x^3 + 546x^2              $&$ 0.011431$ \\
$y^2 + 11xy + 1050y = x^3 + 210x^2                $&$ 0.011556$ \\
$y^2 - 65xy - 1518y = x^3 - 1518x^2               $&$ 0.011851$ \\
$y^2 - 156xy - 268800y = x^3 + 1920x^2            $&$ 0.012069$ \\
$y^2 - 17xy + 960y = x^3 - 30x^2                  $&$ 0.012253$ \\
$y^2 - 45xy + 51030y = x^3 - 810x^2               $&$ 0.012847$ \\
$y^2 - 3141xy - 1765062090y = x^3 + 614790x^2     $&$ 0.013138$ \\
$y^2 - 53xy - 10080y = x^3 + 210x^2               $&$ 0.013244$ \\
$y^2 - 28xy - 840y = x^3 - 420x^2                 $&$ 0.013417$ \\
$y^2 + 535xy + 10929600y = x^3 + 22770x^2         $&$ 0.013564$ \\
$y^2 + 99xy + 1928934y = x^3 - 30618x^2           $&$ 0.013645$ \\
$y^2 - 23xy - 1848y = x^3 + 924x^2                $&$ 0.013835$ \\
$y^2 - 65xy - 549450y = x^3 - 12210x^2            $&$ 0.014833$ \\
$y^2 + 211xy + 537030y = x^3 + 6630x^2            $&$ 0.015087$ \\
$y^2 - 133xy - 35190y = x^3 - 7038x^2             $&$ 0.015425$ \\
$y^2 + 61xy - 840y = x^3 - 840x^2                 $&$ 0.015765$ \\
$y^2 - 6525xy + 199905468750y = x^3 - 22781250x^2 $&$ 0.015798$ \\
$y^2 - 80xy - 34560y = x^3 + 480x^2               $&$ 0.015833$ \\
$y^2 - 28xy - 89100y = x^3 - 3300x^2              $&$ 0.016252$ \\
$y^2 - 157xy + 288990y = x^3 - 1170x^2            $&$ 0.016281$ \\
$y^2 + 103xy + 26730y = x^3 + 330x^2              $&$ 0.016335$ \\
$y^2 + 37xy + 2184y = x^3 + 168x^2                $&$ 0.016637$ \\
$y^2 + 19xy + 270y = x^3 + 30x^2                  $&$ 0.016703$ \\
$y^2 + 103xy + 84150y = x^3 + 990x^2              $&$ 0.016966$ \\
$y^2 - 41xy - 630y = x^3 - 630x^2                 $&$ 0.017003$ \\
$y^2 - 103xy - 2522520y = x^3 - 32760x^2          $&$ 0.017111$ \\
$y^2 + 43xy - 5670y = x^3 + 210x^2                $&$ 0.017276$ \\
$y^2 - 21xy + 10890y = x^3 - 990x^2               $&$ 0.017371$ \\
$y^2 - 181xy - 436800y = x^3 + 2730x^2            $&$ 0.017812$ \\
$y^2 + 100xy + 70200y = x^3 + 780x^2              $&$ 0.018474$ \\
$y^2 + 77xy - 220500y = x^3 + 1260x^2             $&$ 0.018503$ \\
$y^2 + 20xy + 53760y = x^3 - 1920x^2              $&$ 0.018833$ \\
$y^2 + 341xy + 2827440y = x^3 + 9240x^2           $&$ 0.018965$ \\
$y^2 + 99xy + 65610y = x^3 + 1458x^2              $&$ 0.019140$ \\
$y^2 - 2xy - 150y = x^3 + 30x^2                   $&$ 0.019248$
\end{tabular}
\end{adjustwidth}

\subsection{Most Integral Points}

After looking at the candidate curves for the previous searches, we noticed that several curves from the $8P$ parameterization had a rank of at least $3$. This is unusual, and it allows us to use the parameterizations to generate curves with many integral points in general. A high rank is desirable because it allows us to add the generators in more combinations without increasing the height as much. 

We performed a new search as follows:
\begin{enumerate}
\item In \emph{Magma}, we searched through the parameterization for $8P$ integral for curves which have a rank of at least $3$ and added them to a new list of candidate curves.
\item Then, in \emph{Magma} we computed the number of integral points for each of these new candidate curves.
\end{enumerate}

The rank of each curve was calculated using \emph{DescentInformation} in \emph{Magma} with \emph{HeightBound} equal to $11$ for each curve. The generators returned were used as the \emph{FBasis} parameter in the \emph{IntegralPoints} routine in \emph{Magma} to obtain the integral points on the curve.

Below are first $20$ curves we found with the most integral points:

\clearpage

\begin{adjustwidth}{-0.75in}{-0.75in}
\begin{tabular}{l|l|l}
Curve                                               & Integral points& Rank\\
$y^2+30031xy+2603255431140y=x^3+115746540x^2       $& $272$ & $5$\\
$y^2-180839xy-78666732144000y=x^3+441451920x^2     $& $244$ & $5$\\
$y^2-151801xy-102037611637500y=x^3+698289900x^2    $& $236$ & $5$\\
$y^2-113051xy-11790474696000y=x^3+105225120x^2     $& $230$ & $5$\\
$y^2-17189xy+2269552398750y=x^3-106177890x^2       $& $204$ & $5$\\
$y^2-1031xy-72558720y=x^3+77520x^2                 $& $202$ & $4$\\
$y^2-23639xy+1912486878720y=x^3-71404080x^2        $& $202$ & $5$\\
$y^2-19289xy-271964385000y=x^3+14922600x^2         $& $194$ & $4$\\
$y^2-21599xy-467606805120y=x^3+23030280x^2         $& $190$ & $4$\\
$y^2-12905xy-100926385560y=x^3+8288280x^2          $& $174$ & $4$\\
$y^2-73601xy-10774570847040y=x^3+151464390x^2      $& $170$ & $4$\\
$y^2-90691xy-35966223556032y=x^3+422327136x^2      $& $170$ & $5$\\
$y^2-74309xy+109358873613990y=x^3-1251260010x^2    $& $168$ & $4$\\
$y^2-5825xy-10446885900y=x^3+1918620x^2            $& $168$ & $4$\\
$y^2-251129xy-211305388694400y=x^3+855127350x^2    $& $168$ & $4$\\
$y^2-9329xy-30288244140y=x^3+5419260x^2            $& $166$ & $4$\\
$y^2-6161xy+14499985500y=x^3-2222220x^2            $& $164$ & $4$\\
$y^2-54791xy-8933331939840y=x^3+175245840x^2       $& $164$ & $4$\\
$y^2-89879xy-38698378266240y=x^3+462545160x^2      $& $164$ & $5$\\
$y^2-302353xy+7331571365922804y=x^3-20602345764x^2 $& $164$ & $5$
\end{tabular}
\end{adjustwidth}

\clearpage

\section {Self-Descriptive Numbers}
\subsection{Introduction}
\label {intro}
One popular puzzle is to create a so-called ``self-descriptive number'' in base 10: a number whose first digit refers to the number of 0's, second digit refers to the number of 1's, third digit the number of 2's, etc. To be more precise, the challenge is to find a 10 digit number $ b_0 b_1 b_2 b_3 b_4 b_5 b_6 b_7 b_8 b_9 $ with $ 0 \leq b_i < 10$ and $b_i = \left | \{b_j | b_j = i\} \right |$.

It is well known that the only solution to this problem in base $10$ is $6210001000$ \cite{keystoinfinity} \cite{wiki}. This problem has also been generalized to larger bases, for which all solutions have been classified. In this paper we will derive the same classification for the general solutions in a new way.
\subsubsection {Problem Description}
The problem considered is, for all possible $n > 0$, to find a number with $n$ digits in base $n$ such that for every $i$ with $0 \leq i < n$, the value of digit $i$ (starting from digit 0) is the number of digits whose value is $i$. The problem can be equivalently formulated without relying on different bases by instead writing each digit as an entry in an ordered list. For example, the hexadecimal number A0B6 could be written as $(10, 0, 11, 6)$, and this is the approach taken when the base is larger than $10$. A self-descriptive number in base $n$ is, for this reason, instead referred to more generally as a \emph{solution} of length $n$. Entries of the solution, however, may still be refered to as digits.
\subsubsection {Notation}
We will frequently refer to n-digit solutions and the digits within them. For indeterminate solutions referenced by a variable such as $c$, the i-th digit is represented using function notation as $c(i)$. For a solution of length $n$ which is not indeterminant, we represent it as a $n$-digit number with digits 0-9 when possible, and we refer to it as a \emph{literal solution}. For instance, in the literal solution 1210, there is one 0, two 1's, one 2, and zero 3's. If a literal solution contains a digit larger than 9, we represent it instead using an ordered list in order to avoid ambiguity. For instance, the following is a literal solution for 15 digits: $$(11, 2, 1, 0, 0, 0, 0, 0, 0, 0, 0, 1, 0, 0, 0)$$ Finally, for any solution $c$, we denote $L(c)$ to be the length of $c$ (or the base of $c$).

\subsection {Bounding Zeros}
Notice in the 10-digit literal solution $6210001000$, there are a lot of 0's. This is no coincidence, and in this section we shall get a handle on why this is the case. To do so, however, we will need to make one important observation about solutions. Suppose we have a solution $b$. With this solution, we can partition the set of digits into the set of all 0's, the set of 1's, the set of 2's, etc. Then, notice that the size of the partition which contains all $i$ digits is $b(i)$ because $b(i)$ represents the number of $i$'s in the number. Thus, if we sum each digit of the number, we must end up with the length of the solution. For example, in the solution we know for 10 digits, $6210001000$, we have: $$6+2+1+0+0+0+1+0+0+0 = 10$$
To summarize, if $b$ is a solution of length $L(b)$, then:
\begin {equation}
\label {eq1}
\sum_{i=0}^{L(b) - 1} {b(i)} = L(b)
\end {equation}
Instead of just adding the digits, there is another helpful way of computing this sum. Notice that in the solution $6210001000$, there are six 0's which don't contribute to the sum. Instead of counting the number of 0's in the solution, we can look at the first digit. Likewise, the second digit tells us that there are two 1's in the solution, each contributing 1 to the sum. And the third digit tells us that there is one digit which contributes 2 towards the sum. And in fact: $$6 \cdot 0 + 2 \cdot 1 + 1 \cdot 2 + 0 \cdot 3 + 0 \cdot 4 + 0 \cdot 5 + 1 \cdot 6 + 0 \cdot 7 + 0 \cdot 8 + 0 \cdot 9 = 10$$
Generalizing this for any list $b$ of length $L(b)$, define $S(b)$ to be the following sum:
\begin {equation}
\label {eq2}
S(b) = \sum_{i = 0}^{L(b) - 1} {b(i) \cdot i}
\end {equation}
Then, if $b$ is a solution, we must have:
\begin {equation}
\label {eq3}
S(b) = L(b)
\end {equation}
Our equation for $S$ now allows us to bound the number of zeros in any solution. Let $b$ be an arbitrary list with $c$ zeros. Then there are $L(b) - c$ digits with a non-zero value. Looking at (\ref {eq2}), the smallest possible value of $S(b)$ occurs when the non-zero digits are at the beginning of the number, since non-zero digits near the end of the number increase the sum more. Additionally, the sum is minimized when all non-zero digits have a value of one. Thus we must have:
$$S(b) \geq \sum_{i = 0}^{L(b) - c - 1} {i} = \frac {(L(b) - c)(L(b) - c - 1)} {2}$$
Now, if $b$ is a solution, then by (\ref {eq3}) we must have:
$$\frac {(L(b) - b(0))(L(b) - b(0) - 1)} {2} \leq L(b)$$
Using the quadratic formula, we find that:
\begin {equation}
\label {eq4}
b(0) \geq L(b) - \frac {1}{2} - \sqrt {2L(b) + \frac {1}{4}}
\end {equation}
This lower bound on the number of zeros in a solution will be the key to proving the uniqueness of solutions of sufficient length.

\subsection {Generating Solutions}
In this section, we will demonstrate how to generate new solutions of different lengths given a particular solution subject to appropriate constraints.
\subsubsection {Infinite Families of Solutions}
\label {sec31}
Let $b$ be some solution. Then we will call $b$ \emph{extendable} if and only if the following conditions hold:
\begin {enumerate}
\item $b(b(0)) = 1$
\item For all $i$ such that $L(b) > i > b(0)$, we have $b(i) = 0$
\end {enumerate}
Additionally, we will call $b$ a \emph{particular solution} if and only if:
\begin {enumerate}
\item $b$ is an extendable solution
\item $b(b(0) - 1) > 0$
\end {enumerate}

\noindent{Example: $72100001000$ Is an extendable solution, but not particular.\\}

If $b$ is not extendable, then we will call $b$ a \emph{sporadic solution}. Extendable solutions were so named because given an extendable solution, we can generate solutions of arbitrarily larger lengths. Suppose $b$ is an extendable solution. Define $e = E(b)$, a list of length $L(b) + 1$, as the \emph{extension} of $b$ by the following:
$$ e(i) = 
\begin {cases}
b(0) + 1,& \text{if } i = 0\\
0,& \text{if } i = b(0)\\
1,& \text{if } i = b(0) + 1\\
0,& \text{if } i = L(b)\\
b(i),& \text{otherwise}
\end {cases}
$$

\noindent{Example: $E(72100001000) = 821000001000$\\}

Importantly, if $b$ is an extendable solution, then $E(b)$ is an extendable solution. Let $b$ be any extendable solution, and $e = E(b)$. Then $e(b(0)) = 0$ and $e(b(0) + 1) = 1$. But since $b$ was extendable, we also have that $b(b(0)) = 1$ and $b(b(0) + 1) = 0$. Since $e(L(b)) = 0$, this means that the number of 0's in $e$ is $b(0) + 1$, and this is the value of $e(0)$. It also means that the number of 1's in $e$ is the same as the number of 1's in $b$, and indeed we have $e(1) = b(1)$. The number of digits with the value $b(0)$ in $b$ is one, so no digits have the value $b(0)$ in $e$. Correspondingly, $e(b(0)) = 0$. Likewise, the number of digits with the value $b(0) + 1$ in $b$ is zero, and exactly one digit has the value $b(0) + 1$ in $e$. Indeed, we have that $e(b(0) + 1) = 1$. For any remaining value $d$, the number of digits with value $d$ in $e$ is the same as the number of digits with value $d$ in $b$, and we also have $e(d) = b(d)$. Thus $e$ is a solution. Further, $e$ is also extendable because we have that $e(e(0)) = e(b(0) + 1) = 1$ and if $i > b(0) + 1$ then $e(i) = b(i) = 0$.

Notice that the extension of a particular solution is not a particular solution because $e(e(0) - 1) = e(b(0)) = 0$. Also, if $e$ is an extendable solution and is not a particular solution, then it is the extension of some other solution $e' = E^{-1}(e)$ defined by:
$$ e'(i) = 
\begin {cases}
e(0) - 1,& \text{if } i = 0\\
0,& \text{if } i = e(0)\\
1,& \text{if } i = e(0) - 1\\
e(i),& \text{otherwise}
\end {cases}
$$
Because of this, we can identify every extendable solution with its unique particular solution which generates it. For instance, the solution $6210001000$ is an extendable (and not particular) solution of length 10, so we can apply $E^{-1}$ until we get a particular solution:
\begin {itemize}
\item $E^{-1}(6210001000) = 521001000$
\item $E^{-1}(521001000) = 42101000$
\item $E^{-1}(42101000) = 3211000$
\end {itemize}
So we have reduced the question of identifying solutions to the question of identifying particular and sporadic solutions.
\subsubsection {Long Solutions}
Now we will show that long solutions must be extendable solutions, and can't be particular.

As an immediate consequence of (\ref {eq4}), for any solution $b$ we have that $b(0) > \lceil\frac {1} {2}L(b)\rceil$ when $L(b) > 11$. Since $b$ is a solution, this means that:
$$b(b(0)) > 0$$
But is it possible that $b(b(0)) > 1$ when $L(b) > 11$? No because if $b(b(0)) > 1$ then by (\ref {eq2}) we would have that:
$$S(b) \geq 2b(0) > 2 \lceil \frac {1} {2}L(b)\rceil \geq L(b)$$
contradicting (\ref {eq3}). And is it possible that $b(b(0) - 1) > 0$? Again, no because this would mean that:
$$S(b) \geq b(0) + b(0) - 1 > \lceil \frac {1} {2}L(b) \rceil + \lceil \frac {1} {2}L(b) \rceil \geq L(b)$$
contradicting (\ref {eq3}). Similarly, if $i > b(0)$ then the digit $b(i)$ must be zero, for if $b(i) > 0$, then by (\ref {eq2}) and the fact that $b(0) > \lceil \frac {1} {2} L(b) \rceil$ we would have:
$$S(b) \geq b(0) + i > \lceil \frac {1} {2}L(b) \rceil + \lceil \frac {1} {2} L(b) \rceil \geq L(b)$$
contradicting (\ref {eq3}). Referring back to section (\ref {sec31}), we have shown that $b$ is an extendable solution, and not a particular solution nor a sporadic solution. This means that it is the extension of some particular solution whose length can be no more than 11.

\subsection {All Solutions}
To finish the classification of solutions, all that is left is to search for particular solutions and sporadic solutions of length less than 11. This is achieved using a python program which searches through all lists satisfying (\ref {eq3}), and finds which ones are solutions. The code is available {\color {blue} \href{https://github.com/been-jamming/Self-Descriptive-Numbers}{here}}. Here are all solutions of length less than 12:

\begin {table}[h!]
\begin {center}

\caption {Small Solutions}
\begin {tabular}{l|c|r}
$L(b)$ & $b$ & Classification\\
\hline
4 & $1210$ & sporadic\\
4 & $2020$ & sporadic\\
5 & $21200$ & sporadic\\
7 & $3211000$ & particular\\
8 & $42101000$ & extendable\\
9 & $521001000$ & extendable\\
10 & $6210001000$ & extendable\\
11 & $72100001000$ & extendable\\
\end {tabular}

\end {center}
\end {table}

Notice that the only particular solution which appears in this table is the one of length 7. This means that every solution of length greater than 11 is the extension of $3211000$. Because extending one solution can only yield one solution of any particular length, this means that for each $n > 11$, there is precisely one solution of the form:
$$n - 4, 2, 1, [n - 7 \text{ zeros...}], 1, 0, 0, 0$$
This concludes the classification of self-descriptive numbers.

\clearpage
\bibliography{full_text}
\bibliographystyle{plain}
\end {document}